\documentclass[a4paper, 11pt]{article}
\usepackage[mathscr]{eucal}
\usepackage{amsmath}
\usepackage{amssymb}
\usepackage{amsfonts}
\usepackage{amsthm}

\usepackage{cases}
\usepackage[mathscr]{eucal}
\usepackage{amsmath}
\usepackage{amssymb}
\usepackage{amsfonts}
\usepackage{amsthm}
\usepackage[dvips]{graphics, color}
\theoremstyle{plain}
\newtheorem{theorem}{Theorem}[section]
\newtheorem{proposition}{Proposition}[section]
\newtheorem{corollary}{Corollary}[section]
\newtheorem{remark}{Remark}[section]

\newtheorem{example}{Example}[section]
\newtheorem{definition}{Definition}[section]
\newtheorem{problem}{Problem}[section]
\newtheorem{conjecture}{Conjecture}[section]

\numberwithin{equation}{section}

\def\<{\left<} \def\>{\right>}

\def\bea{\begin{eqnarray} }
\def\eea{\end{eqnarray} }
\def\be{\begin{equation} }
\def\ee{\end{equation} }
\def\qed{\ifhmode\unskip\nobreak\fi\ifmmode\ifinner\else\hskip5pt
\fi\fi\hbox{\hskip5 pt \vrule width4 pt height6 pt depth1.5 pt \hskip1pt }}

\makeatletter

\renewcommand\section{\@startsection {section}{1}{\z@}%
            {-4.5ex \@plus -1ex \@minus -.2ex}%
            {1.5ex \@plus.2ex}%
            {\normalfont\large\bfseries}}

\renewcommand\subsection{\@startsection {subsection}{1}{\z@}%
            {-3.5ex \@plus -1ex \@minus -.2ex}%
            {1.5ex \@plus.2ex}%
            {\normalfont\normalsize\bfseries}}

\begin{document}

\title{\bf\Large A short survey on $\delta$-ideal CR submanifolds}
\author{\bf Toru Sasahara\\
\small \it Center for Liberal Arts and Sciences,  Hachinohe Institute of Technology\\
\small \it Hachinohe, Aomori 031-8501, JAPAN\\
\small {\it E-mail}: {\tt sasahara@hi-tech.ac.jp}}

\date{}
\maketitle

\begin{abstract}
{\footnotesize This paper surveys some of the known  results  on $\delta$-ideal 
CR submanifolds in complex space forms, the nearly K\"{a}hler $6$-sphere and odd dimensional
unit spheres.  
In addition, 
the relationship between $\delta$-ideal CR submanifolds and critical points of 
the $\lambda$-bienergy is mentioned.
Some topics on variational problem for the $\lambda$-bienergy 
 are also presented.
}
\end{abstract}




 \section{Introduction}

In submanifold theory,  it is important to establish relations between 
  extrinsic and intrinsic invariants of submanifolds.
   %
In the early 1990s, the notion of  $\delta$-invariants was introduced by Chen
(see \cite{chen2}, \cite{chen3} and \cite{chen4}).
These invariants are obtained by subtracting a certain amount of sectional curvatures 
from the scalar curvature.
Furthermore, 
 he established pointwise optimal 
 inequalities involving
 $\delta$-invariants and 
 the squared mean curvature
  of 
 arbitrary submanifolds in real and complex space forms.
A submanifold is   said to be {\it $\delta$-ideal} if it it satisfies an
 equality case of the inequalities everywhere. 
During the last two decades,  
many interesting results on $\delta$-ideal submanifolds have been obtained.
%

The main purpose of  this paper is to survey
some of the known  results  on $\delta$-ideal CR submanifolds in 
  complex space forms, the nearly K\"{a}hler $6$-sphere and odd dimensional unit spheres. 
 For a given compact almost CR manifold $M$ (with or without boundary)
 equipped with 
 a compatible metric, the $\delta$-ideal CR 
 immersions of $M$ minimize the  $\lambda$-bienergy 
 among all isometric CR immersions of $M$.  
In view of this fact,  
some topics on variational problem for the $\lambda$-bienergy are also
  presented.

\section{Preliminaries}
Let $M$ be  
 an $n$-dimensional  submanifold of a Riemannian manifold $\tilde M$.
 Let us denote by $\nabla$ and
 $\tilde\nabla$ the Levi-Civita connections on $M$ and 
 $\tilde M$, respectively. The
 Gauss and Weingarten formulas are respectively given by
\be
 \begin{split}
 \tilde \nabla_XY&= \nabla_XY+B(X,Y), \label{gawe}\\
 \tilde\nabla_X V&= -A_{V}X+D_XV
 \end{split}\nonumber
\ee
 for tangent vector fields $X$, $Y$ and normal vector field $V$,
 where $B$, $A$ and $D$ are the second fundamental 
 form, the shape operator and the normal
 connection.

 The mean curvature vector field $H$ is defined by 
$H=(1/n){\rm trace}\hskip2pt B.$
 The function $|H|$ is called the  {\it mean curvature}.
 If it vanishes identically, then $M$ is called a {\it minimal submanifold}.
 In particular, if $B$ vanishes identically,  then $M$ is called a {\it totally geodesic submanifold}.

 
\begin{definition}[\cite{be}]
{\rm   Let $M$ be a Riemannian submanifold  of an almost Hermitian  manifold
$\tilde M$ and let $J$ be the complex  structure of $\tilde M$. A submanifold $M$
 is called a {\it CR submanifold} if there exist differentiable distributions
 $\mathcal{H}$ and  ${\mathcal H}^{\perp}$ such that 
 \be 
 TM=\mathcal{H}\oplus {\mathcal H}^{\perp}, \quad \nonumber
 J\mathcal{H}=\mathcal{H}, \quad
 J{\mathcal H}^{\perp}\subset T^{\perp}M,
 \ee 
where $T^{\perp}M$
denotes the normal bundle of $M$. 
 A CR submanifold is called a {\it K\"{a}hler submanifold} (resp. {\it totally real submanifold})
  if {\rm rank}\hskip2pt${\mathcal H}^{\perp}=0$ (resp. {\rm rank}\hskip2pt${\mathcal H}=0$).
A totally real submanifold is called a {\it Lagraingian submanifold} if 
$J(TM)=T^{\perp}M$.
A CR submanifold is said to be ${\it proper}$ if {\rm rank}\hskip2pt${\mathcal H}\ne 0$ 
and {\rm rank}\hskip2pt${\mathcal H}^{\perp}\ne 0$.
}
\end{definition}

\section{$\delta$-invariants}
Let $M$ be an $n$-dimensional Riemannian manifold. Denote by $K(\pi)$ the sectional curvature
 of $M$ associated with a plane section $\pi\subset T_pM$, $p\in M$. For any orthonormal basis
 $\{e_1, \ldots, e_n\}$ of the tangent space $T_pM$, the scalar curvature $\tau$ at $p$ is defined
 by
 \be
 \tau(p)=\sum_{i<j}K(e_i \wedge e_j).\nonumber
 \ee
 Let $L$ be a subset of $T_pM$ of dimension $r\geq 2$ and $\{e_1, \ldots, e_r\}$
 an orthonormal basis of $L$. We define the scalar curvature $\tau(L)$ of the $r$-plane section $L$ by
 \be
 \tau(L)=\sum_{\alpha<\beta}
K(e_{\alpha}\wedge e_{\beta}), \nonumber
\quad 1\leq\alpha, 
\beta\leq r.
 \ee
 For an integer $k\geq 0$, denote by $\mathcal{S}(n, k)$ the finite set which consists of unordered
$k$-tuples $(n_1, \ldots, n_k)$ of integers satisfying $2\leq n_1 \ldots, n_k<n$ and 
$n_1+\ldots+n_k\leq n$.
 We denote by $\mathcal{S}(n)$ the set of $k$-tuples with $k\geq 0$ for a fixed $n$.
 
For each $k$-tuple $(n_1, \ldots, n_k)\in \mathcal{S}(n)$, 
  the notion
 of {\it $\delta$-invariant} $\delta(n_1, \ldots, n_k)$ was introduced by  Chen \cite{chen3}
  as follows:
\be
\delta(n_1, \ldots, n_k)(p)=\tau(p)-{\rm inf}\{\tau(L_1)+\cdots+\tau(L_k)\},\nonumber
\ee
where $L_1$, \ldots, $L_k$ run over all $k$ mutually orthogonal subspaces of $T_pM$ such that
 $\dim L_j=n_j$, $j=1, \ldots, k$.

Let  $\overline{Ric}$ denote the maximum Ricci curvature function on $M$ defined by
\be\overline{Ric}(p)={\rm max}\{S(X, X)|X\in U_pM\},\nonumber\ee
where $S$ is the Ricci tensor and $U_pM$
is the unit tangent vector space of $M$ at $p$.
Then, we have $\delta(n-1)(p)=\overline{Ric}(p)$.
 
 
Let $M$ be a K\"{a}hler manifold with real dimension $2n$.
 For each $k$-tuple $(2n_1, \ldots, 2n_k)\in \mathcal{S}(2n)$,  Chen \cite{chen3} also 
 introduced the notion of 
  {\it complex}
$\delta$-$invariant$ $\delta^c(2n_1, \ldots, 2n_k)$, which is defined by
\be
\delta^c(2n_1, \ldots, 2n_k)(p)=\tau(p)-{\rm inf}\{\tau(L_1^c)+\cdots+\tau(L_k^c)\},\nonumber
\ee
where $L_1^c$, \ldots, $L_k^c$ run over all $k$ mutually orthogonal complex
subspaces of $T_pM$ such that
 ${\rm dim}L_j=2n_j$, $j=1, \ldots, k$.
 
For simplicity, we denote $\delta(\lambda, \ldots, \lambda)$ and 
$\delta^c(\lambda, \ldots, \lambda)$  
 by $\delta_k(\lambda)$ and $\delta_k^c(\lambda)$, 
respectively, where $\lambda$ appears $k$ times.


 \section{Inequalities involving $\delta$-invariants and ideal  submanifolds}

 
 For each $(n_1, \ldots, n_k)\in{\mathcal S}(n)$,  let $c(n_1, \ldots, n_k)$ and 
 $b(n_1, \ldots, n_k)$ be the constants given by
\begin{align}
&c(n_1, \ldots, n_k)=\frac{n^2(n+k-1-\sum_{j=1}^k n_j)}{2(n+k
-\sum_{j=1}^k n_j)},\nonumber\\
&b(n_1, \ldots, n_k)=\frac{1}{2}\Bigl(n(n-1)-\sum_{j=1}^k n_j(n_j-1)\Bigr).\nonumber
\end{align}
 
Chen obtained the following inequality
 for  an arbitrary submanifold in a
  real space form.
 \begin{theorem}[\cite{chen4}]\label{ingere}
Given an $n$-dimensional submanifold $M$ in  an $m$-dimensional  real space form
$R^m(\epsilon)$ of constant sectional curvature $\epsilon$, we have
 \be
 \delta(n_1, \ldots, n_k)
 \leq c(n_1, \ldots, n_k)|H|^2
+b(n_1, \ldots, n_k)\epsilon.\label{inre}
 \ee
 Equality sign of $(\ref{inre})$ holds at a point $p\in M$ for some $(n_1, \ldots, n_k)\in{\mathcal S}(n)$ 
if and only if there exists an orthonormal basis $\{e_1, \ldots, e_{2m}\}$ at $p$ 
such that $e_1, \ldots, e_{n}$ are tangent to $M$ and  the shape operators of $M$ in $R^m(\epsilon)$ at $p$ take  following forms$:$
\be
A_{e_r}= \left(
    \begin{array}{cccc}
      A_1^r & \ldots & 0 & \\
      \vdots & \ddots & \vdots &  ${\rm \mbox{\Large 0}}$\\
      0 & \ldots & A_k^r &  \\
       &   ${\rm \mbox{\Large 0}}$ &  & \mu_r  I
    \end{array}
  \right),\label{A}
\ee
$$r=n+1, \ldots, 2m,$$
where each $A_{j}^r$ is a symmetric $n_j\times n_j$ submatrix such that
\be
{\rm trace}(A_{1}^r)=\cdots={\rm trace}(A_{k}^r)=\mu_r.\label{traceA}
\ee
 \end{theorem}

 Let $\tilde M^m(4\epsilon)$ be a complex space form of complex dimension $m$ and constant
 holomorphic sectional curvature $4\epsilon$ and let $J$ be the complex
structure of $\tilde M^m(4\epsilon)$.

Let $M$ be an $n$-dimensional submanifold in $\tilde M^m(4\epsilon)$.
 For any vector $X$ tangent to $M$, we put $JX=PX+FX$, where $PX$ and $FX$ are tangential 
 and normal components of $JX$, respectively. For a subspace $L\subset T_pM$ of dimension $r$,
 we set
\be
\Psi(L)=\sum_{1\leq i<j\leq r}\<Pu_i, u_j\>^2,\nonumber 
\ee
 where $\{u_1, \ldots, u_r\}$ is an orthonormal basis of $L$.
 
 For an arbitrary  submanifold in a complex space form, we have
 \begin{proposition}[\cite{chen4}]\label{pro1}
 Let $M$ be 
 an $n$-dimensional submanifold  in a complex space form $\tilde M^m(4\epsilon)$.
 Then, for mutually orthogonal subspaces $L_1$, \ldots, $L_k$ of $T_pM$ such that
 $\dim L_j=n_j$,  we have
 \be
 \tau-\sum_{i=1}^k\tau(L_i)\leq c(n_1, \ldots, n_k)|H|^2+b(n_1, \ldots, n_k)\epsilon
+\frac{3}{2}|P|^2\epsilon
 -3\epsilon\sum_{i=1}^k\Psi(L_i).\label{gein}
 \ee
The equality case of inequality $(\ref{gein})$
 holds at a point $p\in M$ if and only if there exists an orthonormal basis $\{e_1, \ldots, e_{2m}\}$
 at $p$ such that

${\rm (a)}$ $L_j={\rm Span}\{e_{n_1+\cdots+n_{j-1}+1}, \ldots,
 e_{n_1+\cdots+n_j}\}$, \quad $j=1, \ldots, k$,

${\rm (b)}$ the shape operators of $M$ in $\tilde M^m(4\epsilon)$ at $p$ satisfy {\rm (\ref{A})} and 
{\rm (\ref{traceA})}.
 \end{proposition}

 Using Proposition \ref{pro1}, we obtain the following  inequalities.
\begin{proposition}[\cite{chen4}]\label{proka}
Let $M$ be a K\"{a}hler submanifold with real dimension $2n$ in a complex space form
 $\tilde M^m(4\epsilon)$. Then,
we have 
 \be
\delta^c(2n_1, \ldots, 2n_k)\leq 2\biggl(n(n+1)-\sum_{j=1}^kn_j(n_j+1)\biggr)\epsilon.\label{ka}
\ee
The equality case of inequality $(\ref{ka})$ holds at a point $p\in M$ if and only if
there exists an orthonormal basis
 $\{e_1, \ldots, e_{2m}\}$ at $p$ such that $e_1, \ldots, e_{2n}$ are tangent to $M$
 and $e_{2l}=Je_{2l-1}$ $(1\leq l\leq k)$, and moreover, the shape operators of $M$ in 
 $\tilde M^m(4\epsilon)$ at $p$ take
 the following forms$:$
\be
A_{e_r}= \left(
    \begin{array}{cccc}
      A_1^r & \ldots & 0 & \\
      \vdots & \ddots & \vdots &  ${\rm \mbox{\Large 0}}$ \\
       0& \ldots & A_k^r &  \\
       & ${\rm \mbox{\Large 0}}$&  & ${\rm \mbox{\Large 0}}$
    \end{array}
  \right),\nonumber
\ee
$$r=2n+1, \ldots, 2m,$$ where each $A_j^r$ is a symmetric $(2n_j)\times(2n_j)$ submatrix satisfying
${\rm trace}(A_j^r)=0$.
\end{proposition}


 \begin{proposition}[\cite{sa3}]\label{pro3}
 Let  $M$ be an $n$-dimensional CR submanifold with ${\rm rank}\hskip2pt{\mathcal H}=2h$
in 
${\mathbb C}H^m(-4)$.
 Then, we have
 \be
 \delta(n_1, \ldots, n_k)
 \leq c(n_1, \ldots, n_k)|H|^2
-b(n_1, \ldots, n_k)-3h+\frac{3}{2}\sum_{j=1}^k n_j.\label{inCH}
 \ee
 Equality sign of $(\ref{inCH})$ holds at a point $p\in M$ for some $(n_1, \ldots, n_k)\in{\mathcal S}(n)$ 
if and only if there exists an orthonormal basis $\{e_1, \ldots, e_{2m}\}$ at $p$ such that

${\rm (a)}$ each $L_j={\rm Span}\{e_{n_1+\cdots+n_{j-1}+1}, \ldots, e_{n_1+\cdots+n_j}\}$
satisfies 
$\Psi(L_j)=n_j/2$ for $1\leq j\leq k$, 

${\rm (b)}$ the shape operators of $M$ in ${\mathbb C}H^m(-4)$ at $p$ satisfy {\rm (\ref{A})} and {\rm (\ref{traceA})}.
 \end{proposition}

\begin{proposition}[\cite{sa4}]\label{pro4}
Let  $M$ be an $n$-dimensional CR submanifold with ${\rm rank}\hskip2pt{\mathcal H}=2h$
in 
${\mathbb C}P^m(4)$.
 Then, we have
 \be
 \delta(n_1, \ldots, n_k)
 \leq c(n_1, \ldots, n_k)|H|^2
+b(n_1, \ldots, n_k)+3h.\label{inCP}
 \ee
 Equality sign of $(\ref{inCP})$ holds at a point $p\in M$ for some $(n_1, \ldots, n_k)\in{\mathcal S}(n)$ 
if and only if there exists an orthonormal basis $\{e_1, \ldots, e_{2m}\}$ at $p$ such that

${\rm (a)}$ each $L_j={\rm Span}\{e_{n_1+\cdots+n_{j-1}+1}, \ldots, e_{n_1+\cdots+n_j}\}$
satisfies
$\Psi(L_j)=0$ for $1\leq i \leq k$,

${\rm (b)}$ the shape operators of $M$ in ${\mathbb C}P^m(4)$ at $p$ satisfy $(\ref{A})$ and
 $(\ref{traceA})$.
\end{proposition}


\begin{definition}{\rm
A submanifold is 
said to be {\it $\delta(n_1, \ldots, n_k)$-ideal} 
if it satisfies the equality case of  (\ref{inre}), (\ref{inCH}) or (\ref{inCP})  identically  for 
a $k$-tuple $(n_1, \ldots, n_k)\in {\mathcal S}(n)$. Similarly, 
a K\"{a}hler submanifold is said to be {\it $\delta^c(2n_1, \ldots, 2n_k)$-ideal} if it
 satisfies the equality case of (\ref{ka}) identically for 
a $k$-tuple $(2n_1, \ldots, 2n_k)\in {\mathcal S}(2n)$.}
\end{definition}

For more information on $\delta$-invariants and $\delta$-ideal submanifolds, 
we refer the reader to \cite{chen}.

\begin{definition}{\rm
A submanifold is said to be {\it linearly full} in $\tilde M^m(4\epsilon)$
if it  does not lie in any totally geodesic K\"ahler hypersurfaces of $\tilde M^m(4\epsilon)$.}
\end{definition}


 \section{Ideal CR submanifolds in complex hyperbolic space}

We first recall some basic definitions on hypersurfaces.
\begin{definition}{\rm Let $N$ be a submanifold in a Riemannian manifold $\tilde M$ and $UN^{\perp}$
 the unit normal bundle of $N$.
Then, for a sufficiently small $r>0$, the following mapping is an  immersion:
\be
f_r: UN^{\perp}\rightarrow \tilde M,  \quad f_r(p, V)={\exp}_p(rV),\nonumber
\ee
where exp denotes the exponential mapping of $\tilde M$. The hypersurface 
$f_r(UN^{\perp})$ of $\tilde M$  is called the {\it tubular hypersurface} over $N$ with radius $r$.
If $N$ is a point $x$ in $\tilde M$, then the tubular hypersurface over $x$  is a geodesic hypersphere centered
at $x$.}
\end{definition}

\begin{definition}
{\rm For a given point $p\in\mathbb{C}H^m(-4)$, let $\gamma(t)$ be a geodesic with $\gamma(0)=p$, which is parametrized 
by arch length. Denote by $S_t(\gamma(t))$ the geodesic hypersphere
centered at $\gamma(t)$ with radius $t$. 
The limit of $S_t(\gamma(t))$ when $t$ tends to infinity is called a {\it  horosphere}.}
 \end{definition}

 
 \begin{definition}{\rm 
 Let $M$ be a
  real hypersurface  in an almost Hermitian manifold and $V$ be a unit normal vector.
  A hypersurface $M$ is called a {\it Hopf hypersurface} 
 if $JV$ is a principal curvature vector.}
 \end{definition}

 A real hypersurface in an almost Hermitian manifold is a proper CR submanifold with 
 ${\rm rank}{\hskip 2pt}{\mathcal H}^{\perp}=1$.
 The following theorem characterizes the horosphere of $\mathbb{C}H^m(-4)$ in terms of $\delta_k(2)$.
\begin{theorem}[\cite{chen4}]\label{horo1}
Let $M$ be a $\delta_k(2)$-ideal 
real hypersurface of $\mathbb{C}H^m(-4)$.
Then $k=m-1$ and $M$ is an open portion of the horosphere in  $\mathbb{C}H^m(-4)$.
\end{theorem} 
 
\begin{remark}
{\rm The third case of  (9.5) in  \cite{chen4} does not occur, because $L_1\ldots, L_k$  are
 complex planes. Therefore, case (1) of Theorem 9.1 in \cite{chen4}
 shall be removed from the list of $\delta_k(2)$-ideal real hypersurfaces in $\mathbb{C}H^m(-4)$.}
\end{remark}

  For $\delta(2m-2)$-ideal real hypersurfaces
  in $\mathbb{C}H^m(-4)$, Chen proved the following.
\begin{theorem}[\cite{chen5}]\label{2m-2}
Let $M$ be a 
real hypersurface of $\mathbb{C}H^m(-4)$.
Then $M$ is $\delta(2m-2)$-ideal if and only if $M$ is a Hopf hypersurface
with constant mean curvature given by $2\alpha/(2m-1)$, where 
$A_VJV=\alpha JV$ for a unit normal vector $V$.
If $M$ has constant principal curvatures, then 
$M$ is an open portion of one of the following real hypersurfaces:

{\rm (1)} the horosphere of $\mathbb{C}H^2(-4)${\rm ;}

{\rm (2)} the tubular hypersurface over totally geodesic  $\mathbb{C}H^{m-1}(-4)$
in $\mathbb{C}H^m(-4)$ with radius $r=\tanh^{-1}(1/\sqrt{2m-3})$, where $m\geq 3$.
\end{theorem} 
 
It was proved in \cite{chen5} that 
if $m=2$ in Theorem \ref{2m-2}, then
 the
 assumption of  the constancy of principal curvatures 
  is satisfied.  That is to say, we have
\begin{corollary}[\cite{chen5}]
Let $M$ be a $\delta(2)$-ideal real hypersurface of $\mathbb{C}H^2(-4)$. Then
$M$ is an open portion of the horosphere.
\end{corollary}

Let ${\mathbb C}_1^{m+1}$ be the complex number $(m+1)$-space endowed with the complex coordinates
$(z_0, \ldots, z_m)$, the
pseudo-Euclidean metric given by $\tilde g=-dz_0d\bar{w_0}+\sum_{i=1}^mdz_id\bar{w_i}$ and the 
standard complex structure.
For $\epsilon<0$, we put
 $H^{2m+1}_1(\epsilon)=\{z\in {\mathbb C}_1^{m+1}|\<z, z\>=1/\epsilon\}$, 
 where $\< , \>$
denotes the inner product on ${\mathbb C}_1^{m+1}$ induced from $\tilde g$.
For a given $z\in H^{2m+1}_1(\epsilon)$, we put
 $[z]=\{\lambda z| \lambda\in\mathbb C, \lambda\bar{\lambda}=1\}$.
The Hopf fibration is given by 
\be \varpi_{\{m, \epsilon\}}: H^{2m+1}_1(\epsilon)\rightarrow {\mathbb C}H^m(4\epsilon):
 z\mapsto [z]. 
\nonumber\ee



For $\delta_k(2)$-ideal proper CR submanifolds   in ${\mathbb C}H^m(-4)$ 
whose codimensions are greater than one, we have the following representation formula.
\begin{theorem}[\cite{sa1}]\label{horo2}
Let $M$ be a linearly full $(2n+1)$-dimensional $\delta_k(2)$-ideal 
CR submanifold in ${\mathbb C}H^m(-4)$ such that {\rm rank}\hskip2pt${\mathcal H}^{\perp}=1$, $k\geq 1$ and 
$m>n+1$. 
Then,  up to holomorphic isometries of ${\mathbb C}H^m(-4)$, 
 the immersion of $M$ into ${\mathbb C}H^m(-4)$ is given by 
 the composition $\varpi_{\{m, -1\}}\circ z$, where
 \be
 z=\biggl(-1-\frac{1}{2}|\Psi|^2+iu, -\frac{1}{2}|\Psi|^2+iu, \Psi\biggr)e^{it},
 \label{CR1}\ee
and $\Psi$ is a  $2n$-dimensional $\delta^c_n(2)$-ideal K\"{a}hler submanifold in ${\mathbb C}^{m-1}$.
\end{theorem}

Up to holomorphic isometries of ${\mathbb C}H^m(-4)$, the horosphere
 in ${\mathbb C}H^m(-4)$ is a real hypersurface defined by $\{[z]: z\in H^{2m+1}_1(-1), 
 \hskip5pt  
  |z_0-z_1|=1\}$ (see, for example,
 \cite{ver}).
 Hence, Theorem \ref{horo2} can be considered as an extension of Theorem
 \ref{horo1}.

As an immediate corollary of Theorem \ref{horo2}, we obtain
\begin{corollary}[\cite{cv}]
Let $M$ be a linearly full $3$-dimensional $\delta(2)$-ideal 
CR submanifold in ${\mathbb C}H^m(-4)$ such that {\rm rank}\hskip2pt${\mathcal H}^{\perp}=1$
and 
$m>2$. 
Then,  up to holomorphic isometries of ${\mathbb C}H^m(-4)$, 
 the immersion of $M$ into ${\mathbb C}H^m(-4)$ is given by
 the composition $\varpi_{\{m, -1\}}\circ z$, where $z$ is given by {\rm (\ref{CR1})} and 
 $\Psi(w)$ 
is a  holomorphic curve in ${\mathbb C}^{m-1}$ with $\Psi^{\prime}(w)\ne 0$.
\end{corollary}

For  general 
$\delta(n_1, \ldots, n_k)$-ideal  proper CR submanifolds   in ${\mathbb C}H^m(-4)$ 
whose codimensions are greater than one, the following classification result has been obtained.
\begin{theorem}[\cite{sa3}, \cite{sa4}]\label{cla}
Let $M$ be a linearly full $(2n+1)$-dimensional $\delta(n_1, \ldots, n_k)$-ideal 
CR submanifold in ${\mathbb C}H^m(-4)$ such that {\rm rank}\hskip2pt${\mathcal H}^{\perp}=1$, $k\geq 1$ and 
$m>n+1$. Then, we have
$JH\in\mathcal{H}^{\perp}$, $A_VJV=(2n/\sqrt{k(2n-k)})JV$ for $V=H/|H|$,  $DH=0$,  and moreover,
  the mean curvature is given by $$\frac{2n(k+1)}{(2n+1)\sqrt{k(2n-k)}}.$$
  If 
all principal curvatures of $M$ 
with respect to $H/|H|$ are constant, then one of  the following two cases occurs{\rm :}

{\rm (1)} 
 $M$ is locally congruent with the immersion described in Theorem {\rm \ref{horo2}}.

{\rm (2)} 
$n/k\in {\mathbb Z}-\{1\}$, 
$n_1=\cdots=n_k=2n/k$, and $M$ is locally congruent with the immersion 
\be
\varpi_{\{m, -1\}}\biggl(\varpi_{\{m-1, \frac{2k-2n}{2n-k}\}}^{-1}(\Psi), \sqrt{\frac{k}{2n-2k}}e^{it}\biggr),\nonumber\ee
where 
$\Psi$ is a $2n$-dimensional $\delta^c_k(2n/k)$-ideal K\"{a}hler submanifold in 
${\mathbb C}H^{m-1}
(\frac{8k-8n}{2n-k})$.
\end{theorem}


If $n>1$, $k=1$ and $n_1=2n$ in Theorem \ref{cla}, then we have
\begin{corollary}[\cite{sa2}, \cite{sa4}]\label{2n}
Let $M$ be a linearly full $(2n+1)$-dimensional CR submanifold 
 in ${\mathbb C}H^m(-4)$ such that {\rm rank}\hskip2pt${\mathcal H}^{\perp}=1$, $n>1$ and $m>n+1$.
 Then $M$ is  $\delta(2n)$-ideal if and only if $JH\in\mathcal{H}^{\perp}$,  
 $DH=0$, $A_VJV=(2n/\sqrt{(2n-1)})JV$ for $V=H/|H|$, and moreover,
  the mean curvature is given by $$\frac{4n}{(2n+1)\sqrt{2n-1}}.$$
  If 
all principal curvatures of $M$ 
with respect to $H/|H|$ are constant,   then, up to holomorphic isometries of ${\mathbb C}H^m(-4)$, the immersion of $M$ into ${\mathbb C}H^m(-4)$
is given by
\be
\varpi_{\{m, -1\}}\Biggl(\varpi_{\{m-1, \frac{2-2n}{2n-1}\}}^{-1}(\Psi), \sqrt{\frac{1}{2n-2}}e^{it}\Biggr),
\nonumber\ee
where 
$\Psi$ is a  $2n$-dimensional K\"{a}hler submanifold in ${\mathbb C}H^{m-1}
(\frac{8-8n}{2n-1})$.
\end{corollary}

A hypersurface given by (2) in Theorem \ref{2m-2} can be  rewritten as follows (see, for example, 
\cite[Example 6.1]{mon}):
 \be
\varpi_{\{m, -1\}}\Biggl(H_1^{2m-1}\biggl(\frac{4-2m}{2m-3}\biggr)\times 
S^1\biggl(\frac{1}{\sqrt{2m-4}}\biggr)\Biggr),
\nonumber\ee
where $S^1(r)=\{z\in{\mathbb C}| z\bar{z}=r^2\}$.
Thus, Corollary \ref{2n} can be regarded as an extension of Theorem \ref{2m-2}. 

\bigskip

Let $N$ be a K\"ahler hypersurface with real dimension $2n$ in a complex space form.
Let $V$ and $JV$ be normal vector fields of $N$.
Since $A_{JV}=JA_{V}$ and $JA_V=-A_VJ$ holds (cf. \cite[p.175]{kn}), there exists an 
orthonormal basis $\{e_1, Je_{1}, \ldots, e_{n}, Je_{n}\}$ of $T_{p}N$ with respect to which
the shape operators $A_{V}$ and $A_{JV}$ take the following forms:
\be
A_{V}= \left(
    \begin{array}{ccccc}
      \lambda_1 &  &  & & ${\rm \mbox{\Large 0}}$\\
      & -\lambda_1 & & & \\
      &  & \ddots &  &  \\
      &  &  & \lambda_n & \\
     ${\rm \mbox{\Large 0}}$  &   &  & & -\lambda_n
    \end{array}
  \right), \quad
A_{JV}= \left(
    \begin{array}{ccccc}
      0& \lambda_1 &  & & ${\rm \mbox{\Large 0}}$\\
      \lambda_1& 0 & & & \\
      &  & \ddots &  &  \\
      &  &  & 0 & \lambda_n\\
     ${\rm \mbox{\Large 0}}$  &   &  & \lambda_n & 0
    \end{array}
  \right).\nonumber
\ee

Hence, it follows from Proposition \ref{proka}  that every K\"ahler hypersurface with
real dimension $2n$
in a complex space form is $\delta^c_k(2n/k)$-ideal for any natural number 
$k$ such that  $n/k\in {\mathbb Z}$.   Accordingly, applying Theorem \ref{horo2}
yields the following.
\begin{corollary}[\cite{sa1}]\label{m=n+2}
Let $M$ be a linearly full $(2n+1)$-dimensional $\delta_k(2)$-ideal 
CR submanifold in ${\mathbb C}H^{n+2}(-4)$ 
such that {\rm rank}\hskip2pt${\mathcal H}^{\perp}=1$ and $k\geq 1$. 
Then,  up to holomorphic isometries of ${\mathbb C}H^{n+2}(-4)$, 
 the immersion of $M$ into ${\mathbb C}H^{n+2}(-4)$ is given by 
 the composition $\varpi_{\{n+2, -1\}}\circ z$, where
 $z$ is given by {\rm (\ref{CR1})}, 
and $\Psi$ is a  K\"{a}hler hypersurface in ${\mathbb C}^{n+1}$.
\end{corollary}

Similarly, we obtain the following corollary of Theorem \ref{cla}.
\begin{corollary}
Let $M$ be a linearly full $(2n+1)$-dimensional $\delta(n_1, \ldots, n_k)$-ideal 
CR submanifold in ${\mathbb C}H^{n+2}(-4)$ such that {\rm rank}\hskip2pt${\mathcal H}^{\perp}=1$ and $k\geq 1$.
 If 
all principal curvatures of $M$ 
with respect to $H/|H|$ are constant, then one of the following two cases occurs{\rm :} 
 

{\rm (1)}
$M$ is locally congruent with the immersion described in Corollary {\rm \ref{m=n+2}}.


{\rm (2)} $n/k\in {\mathbb Z}-\{1\}$, 
$n_1=\cdots=n_k=2n/k$, and $M$ is locally congruent with the immersion 
\be
\varpi_{\{n+2, -1\}}\biggl(\varpi_{\{n+1, \frac{2k-2n}{2n-k}\}}^{-1}(\Psi), \sqrt{\frac{k}{2n-2k}}e^{it}\biggr),\nonumber\ee
where 
$\Psi$ is a  K\"{a}hler hypersurface   in 
${\mathbb C}H^{n+1}
(\frac{8k-8n}{2n-k})$.
\end{corollary}


It is natural to ask the following problem.
\begin{problem}
Find $\delta(n_1, \ldots, n_k)$-ideal 
CR submanifolds with {\rm rank}\hskip2pt${\mathcal H}^{\perp}=1$
in $\mathbb{C}H^m(-4)$ such that the  principal curvatures 
with respect to $H/|H|$ are not all constant.
\end{problem}

Generally, 
$\delta(n_1, \ldots, n_k)$-ideal proper CR submanifolds
 in ${\mathbb C}H^m(-4)$ have the following properties.
\begin{theorem}[\cite{sa3}]
Let $M$ be a linearly full $(2n+q)$-dimensional $\delta(n_1, \ldots, n_k)$-ideal 
CR submanifold in ${\mathbb C}H^m(-4)$ such that ${\rm rank}\hskip2pt{\mathcal H}^{\perp}=q$.
If $q>1$, then $M$ is minimal. If $q=1$ and $m>n+1$, then $M$ is non-minimal and 
satisfies $DH=0$.
\end{theorem}

A differentiable  manifold $M$ is called an {\it almost contact manifold}  if it
admits a unit vector field $\xi$, a one-form $\eta$ and a $(1, 1)$-tensor field $\phi$ satisfying
\be
\eta(\xi)=1, \quad \phi^2=-I+\eta\otimes\xi.\nonumber
\ee
Every almost contact manifold admits a Riemannian metric $g$ satisfying
\be
g(\phi X, \phi Y)=g(X, Y)-\eta(X)\eta(Y).\nonumber
\ee
The quadruplet $(\phi, \xi, \eta, g)$ is called an {\it almost contact metric structure}.

An almost contact metric structure is called a {\it contact metric structure}
 if it satisfies
\be
d\eta(X, Y)=\frac{1}{2}\Bigl(X(\eta(Y))-Y(\eta(X))-\eta([X, Y])\Bigr)=g(X, \phi Y).\nonumber
\ee
A contact metric structure is said to be {\it Sasakian} 
if the tensor field $S$ defined by
\be
S(X, Y)=\phi^2[X, Y]+[\phi X, \phi Y]-\phi[\phi X, Y]-\phi[X, \phi Y]+2d\eta(X, Y)\xi\nonumber
\ee
vanishes identically.
A manifold equipped with a Sasakian structure is called a {\it Sasakian manifold}.
We refer the reader to \cite{boy} for more information on Sasakian manifolds.

Let $M$ be a CR submanifold with {\rm rank}\hskip2pt$\mathcal{H}^{\perp}$
$=1$ in a complex space form. We 
define a one-form $\eta$ by
$\eta(X)=g(U, X)$, where $U$ is a unit tangent vector field lying in 
 $\mathcal{H}^\perp$, and $g$ is an induced metric on $M$. 
We put $\bar U=(1/\sqrt{r})U$, $\bar\eta=\sqrt{r}\eta$ and $\bar g=rg$ for a 
positive constant $r$. Then, the quadruplet 
$(P, \bar U, \bar\eta, \bar g)$ defines an almost  contact structure on $M$ (cf. \cite[p.96]{dj}).

For the almost contact structure $(P, \bar U, \bar\eta, \bar g)$
on a 
CR submanifold described in Theorem {\rm \ref{cla}}, 
we have the following.
\begin{proposition}[\cite{sa4}]
An almost contact structure $(P, \bar U, \bar\eta, \bar g)$ with $r=\sqrt{\frac{k}{2n-k}}$
on a CR submanifold in Theorem {\rm \ref{cla}}. 
becomes a Sasakian structure. In particular, in the case of $(1)$, the structure is Sasakian with respect
to the induced metric.
\end{proposition}




\section{Ideal CR submanifolds in complex projective space}
 All $\delta_k(2)$-ideal 
Hopf hypersurfaces of $\mathbb{C}P^m(4)$ have been determined as follows:
\begin{theorem}[\cite{chen4}]\label{claidCP}
Let $M$ be a $\delta_k(2)$-ideal 
Hopf hypersurface of $\mathbb{C}P^m(4)$.
Then, one of the following three cases occurs{\rm :}

{\rm (1)} $k=1$ and $M$ is an open portion of a geodesic sphere with radius $\pi/4${\rm ;}

{\rm (2)} $m$ is odd, $k=m-1$, and $M$ is an open portion of a tubular hypersurface
with radius $r\in (0, \pi/2)$ over a totally geodesic $\mathbb{C}P^{(m-2)/2}(4)${\rm ;}

{\rm (3)} $m=2$, $k=1$, and $M$ is an open portion of a tubular hypersurface over the
complex quadric curve $Q_1:=\{[z_0,  z_1, z_2]\in \mathbb{C}P^2: z_0^2+z_1^2+z_2^2=0\}$,
 with radius $r=\tan^{-1}((1+\sqrt{5}-\sqrt{2+2\sqrt{5}})/2)=0.33311971\cdots$. Here, 
 $[z_0,  z_1, z_2]$ is a homogeneous coordinate of  $\mathbb{C}P^2$.
\end{theorem}

A real hypersurface  of $\mathbb{C}P^m(4)$ is called a {\it ruled real hypersurface}
if $\mathcal{H}$ is integrable and each leaf of its maximal integral manifolds is 
locally congruent to
$\mathbb{C}P^{m-1}(4)$.
For a unit normal vector $V$ of a ruled real
hypersurface $M$, the shape operator $A_V$ satisfies
\be
A_VJV=\mu JV+\nu U \hskip5pt (\nu\ne 0), \quad
AU=\nu JV, \quad AX=0 \label{ruled}
\ee
for all $X$ orthogonal to both $JV$ and $U$, where
$U$ is a unit vector orthogonal to $JV$, and $\mu$ and $\nu$
are smooth functions on $M$. Thus, all ruled real hypersurfaces of $\mathbb{C}P^m(4)$ are non-Hopf (see \cite{kim2}).

Using Proposition \ref{pro4} and (\ref{ruled}), we find that
every minimal ruled real hypersurface in  $\mathbb{C}P^m(4)$ is
$\delta_k(2)$-ideal for $1\leq k\leq m-1$. Such a hypersurface can be represented as follows:
\begin{theorem}[\cite{ada}]
A minimal ruled hypersurface of  $\mathbb{C}P^m(4)$ is  congruent to
 $\varpi\circ z$, where  $\varpi:S^{2m+1}(1)\rightarrow \mathbb{C}P^m(4)$
is the Hopf fibration and
\be
z(s, t, \theta, w)=e^{\sqrt{-1}\theta}\Bigl(\cos{s}\cos{t}, \cos{s}\sin{t}, (\sin{s})w\Bigr)\nonumber\ee
for $w\in\mathbb{C}^{m-1}$,  $|w|^2=1$, 
$-\pi/2<s<\pi/2$,  $0\leq t,  \theta<2\pi$.
\end{theorem}

It seems interesting to consider the following problem. 
\begin{problem}
Classify $\delta(n_1, \ldots, n_k)$-ideal non-Hopf 
real hypersurfaces in $\mathbb{C}P^m(4)$.
\end{problem}

Let $M$ be an $n$-dimensional $\delta(n_1, \ldots, n_k)$-ideal CR submanifold 
in ${\mathbb C}P^m(4)$. Let $L_j$ be subspaces of $T_pM$ defined in (a) of Proposition \ref{pro4}. 
Define the subspace $L_{k+1}$ by
$L_{k+1}={\rm Span}\{e_{n_1+\cdots+n_k+1}, \ldots, e_{n}\}$.
It is clear that  $T_pM=L_1\oplus\cdots\oplus L_{k+1}$. We denote by $\mathcal{L}_i$ the 
distribution which is   generated by $L_i$. 
Then, we have the following codimension reduction theorem.
\begin{theorem}[\cite{sa4}]\label{core}
Let $M$ be an $n$-dimensional $\delta(n_1, \ldots, n_k)$-ideal CR submanifold with ${\rm rank}\hskip2pt\mathcal{H}^{\perp}=1$
in ${\mathbb C}P^m(4)$.
If $\mathcal{H}^{\perp}\subset \mathcal{L}_i$ for some $i\in\{1, \ldots, k+1\}$, then
$M$ is contained in a totally geodesic K\"{a}hler submanifold ${\mathbb C}P^{\frac{n+1}{2}}(4)$ in  ${\mathbb C}P^m(4)$.
\end{theorem} 

It was proved in \cite{sa4} that
 if $\dim M=3$, then 
the assumption on $\mathcal{H}^{\perp}$  in Theorem \ref{core}
holds. That is  to say, we have the following.
\begin{corollary}[\cite{sa4}]
Let $M$ be a $3$-dimensional $\delta(2)$-ideal proper CR submanifold 
in ${\mathbb C}P^m(4)$.
Then,
$M$ is contained in  ${\mathbb C}P^{2}(4)$.
\end{corollary}

The following problem arises naturally.

\begin{problem}
Find  $\delta(n_1, \ldots, n_k)$-ideal 
CR submanifolds with
 {\rm rank}\hskip2pt${\mathcal H}^{\perp}=1$ 
 in  $\mathbb{C}P^m(4)$ such that the codimensions are greater than one.
\end{problem}

\section{Ideal CR submanifolds in the nearly K\"{a}hler $6$-sphere}
Let $\mathcal O$ be  the Cayley algebra, and denote  by ${\rm Im}\hskip 2pt\mathcal O$
the purely imaginary part of 
 $\mathcal O$.
 We identify ${\rm Im}\hskip 2pt\mathcal O$
  with $\mathbb R^7$ and define the exterior product $u\times v$ on it by
 \be
 u\times v=\frac{1}{2}(uv-vu).\nonumber
 \ee
 The canonical inner product on $\mathbb R^7$
 is given by $\<u ,v\>=-(uv+vu)/2$.

We define the tensor field $J$ of type $(1, 1)$ on $S^6(1)=
\{p\in{\rm Im}\hskip 2pt\mathcal O|
\<, \>=1\}$ by
\be
JX=p\times X \nonumber
\ee
for any $p\in S^6(1)$, $X\in T_pS^6(1)$.  Let $g$ be the standard metric on $S^6(1)$.
Then $(S^6(1), J, g)$ is a
nearly K\"{a}hler manifold, i.e., an almost Hermitian manifold satisfying
$(\nabla_XJ)X=0$ for any $X\in TS^6(1)$, where
$\nabla$ is the Levi-Civita connection with respect to $g$ (cf. \cite[pp.139-140]{kn}).

For  $3$-dimensional $\delta(2)$-ideal proper
CR submanifolds  in the nearly K\"{a}hler $S^6(1)$, we have the following result.
\begin{theorem}[\cite{dv}, \cite{dv2}]
Let $M$ be a $3$-dimensional $\delta(2)$-ideal proper
 CR submanifold in the nearly K\"{a}hler $S^6(1)$.
Then, $M$ is minimal and locally congruent with  the following immersion{\rm :}
\be
\begin{split}
f(t, u, v)=(\cos t\cos u\cos v, \sin t, \cos t\sin u\cos v, \\
\cos t\cos u\sin v, 0, -\cos t\sin u\sin v, 0).\label{3cr}
\end{split}
\ee
\end{theorem}

\begin{remark}{\rm A CR submanifold (\ref{3cr}) can be rewritten as 
\be
x_1^2+x_2^2+x_3^2+x_4^2+x_6^2=1, \quad x_5=x_7=0, \quad 
x_3x_4+x_1x_6=0,\nonumber
\ee
which implies that  it lies in $S^4(1)$.}
\end{remark}

The following theorem determines 
$4$-dimensional  $\delta(2)$-ideal 
proper CR submanifolds in the nearly K\"{a}hler $S^6(1)$.
\begin{theorem}[\cite{ant}, \cite{adv}]
Let $M$ be a $4$-dimensional $\delta(2)$-ideal 
proper CR submanifold in the nearly K\"{a}hler $S^6(1)$.
Then, $M$ is minimal and locally congruent with  the following immersion{\rm :}
\be
\begin{split}
f(t, u, v, w)=(\cos w\cos t\cos u\cos v, \sin w\sin t\cos u\cos v,\\
\sin 2w\sin v\cos u+\cos 2w\sin u, 0, \sin w\cos t\cos u\cos v,\\
\cos w\sin t\cos u\cos v, \cos 2w\sin v\cos u-\sin 2w\sin u).\label{4cr}
\end{split}
\ee
\end{theorem}

\begin{remark}
{\rm A CR submanifold given by (\ref{4cr}) lies in $S^5(1)$.}
\end{remark}

\begin{definition}{\rm 
A $2$-dimensional submanifold $N$
 of the nearly K\"{a}hler $S^6(1)$ is called an {\it almost complex curve}
if $J(T_pN)=T_pN$ for any $p\in N$.}
\end{definition}

Chen has classified $\delta_k(n_1, \ldots, n_k)$-ideal 
Hopf hypersurfaces of  the nearly K\"{a}hler $S^6(1)$ as follows:
\begin{theorem}[{\cite[p.415]{chen}}]\label{hopfs6}
A Hopf hypersurface of  the nearly K\"{a}hler $S^6(1)$ is  $\delta(n_1, \ldots, n_k)$-ideal 
if and only if  it is either

{\rm (1)} a totally geodesic hypersurface, or 

{\rm (2)} an open part of a tubular hypersurface with radius $\pi/2$ over a non-totally geodesic 
almost complex curve of $S^6(1)$. 
\end{theorem}

\begin{remark}
{\rm A tubular hypersurface described in (2) of  Theorem \ref{hopfs6} is a minimal $\delta(\lambda)$-ideal hypersurface for $\lambda\in\{2, 3, 4\}$.}
\end{remark}

It is natural to consider the following problem.
\begin{problem}
Classify $4$-dimensional $\delta(2, 2)$-ideal and  $\delta(3)$-ideal 
CR submanifolds    in the nearly K\"{a}hler $S^6(1)$.
\end{problem}

\section{Ideal contact CR submanifolds in  odd dimensional unit spheres}

For any point $x\in S^{2n+1}(1)\subset \mathbb{C}^{n+1}$, 
we set $\xi=Jx$, where $J$  denotes the canonical
complex structure of  $\mathbb{C}^{n+1}$.
Let $g$ be the standard metric on $S^{2n+1}(1)$ and 
$\eta$ be the one-form given by $\eta(X)=g(X, \xi)$. 
We consider the orthogonal projection $P:T_x\mathbb{C}^{n+1}\rightarrow
T_xS^{2n+1}(1)$.
We define 
a $(1, 1)$-tensor field $\phi$  on $S^{2n+1}(1)$ by 
$\phi=P\circ J$.
Then, the quadruplet $(\phi, \xi, \eta, g)$ is a Sasakian structure (see, for example, \cite{blair}).

\begin{definition}[\cite{yano}]
{\rm   Let $M$ be a Riemannian submanifold  tangent to $\xi$ 
of a Sasakian manifold. A submanifold $M$ 
 is called a {\it contact CR submanifold} if  there exist
 differentiable distributions $\mathcal{H}$ and  ${\mathcal H}^{\perp}$ such that 
 \be 
 TM=\mathbb{R}\xi\oplus\mathcal{H}\oplus {\mathcal H}^{\perp}, \quad \nonumber
 \phi\mathcal{H}=\mathcal{H}, \quad
 \phi{\mathcal H}^{\perp}\subset T^{\perp}M.
 \ee
A contact CR submanifold is said to be ${\it proper}$ if {\rm rank}\hskip2pt${\mathcal H}\ne 0$ 
and {\rm rank}\hskip2pt${\mathcal H}^{\perp}\ne 0$. 
}
\end{definition}

Non-minimal $\delta(2)$-ideal submanifolds in a sphere have been completely described in \cite{da}.
For minimal $\delta(2)$-ideal proper contact CR submanifolds
in $S^{2m+1}(1)$, we have the following codimension reduction theorem.
\begin{theorem}[\cite{mun}]
Let $M^n$ be a minimal $\delta(2)$-ideal proper contact CR submanifold
in $S^{2m+1}(1)$. Then $n$ is even and there exits a totally geodesic Sasakian
$S^{2n+1}(1)$ in  $S^{2m+1}(1)$ containing $M^n$ as a hypersurface.
\end{theorem}

Therefore, it is sufficient to investigate the case of hypersurfaces. 
Let $N$ be a minimal surface in $S^{n}(1)$ and let $UN^{\perp}$ be its unit normal bundle.
Then, a map
\be F: UN^{\perp}\rightarrow S^n(1): V_p\mapsto V_p \nonumber
\ee
is a minimal $\delta(2)$-ideal codimension one immersion (see \cite[Example 9.8]{chen3}).

Munteanu and Vrancken proved the following.

\begin{theorem}[\cite{mun}]
Let $M^{2n}$ be a  minimal $\delta(2)$-ideal proper contact CR hypersurface
in $S^{2n+1}(1)$. Then $M^{2n}$ can be locally considered as 
the unit normal bundle of the Clifford torus
$S^{1}(1/\sqrt{2})\times S^{1}(1/\sqrt{2})\subset S^3(1)\subset S^{2n+1}(1)$. 
\end{theorem}

\section{Related topics}
This section gives an account of the relationship between $\delta(n_1, \ldots, n_k)$-ideal immersions and 
critical points  of  the $\lambda$-bienergy functional $E_{2,\lambda}$. Some topics 
about
 variational problems for $E_{2, \lambda}$  are also presented.
\subsection{$\lambda$-bienergy functional}
Let $f:M \rightarrow N$ be a smooth map of an
 $n$-dimensional Riemannian manifold
into another Riemannian manifold.
The {\it tension field}
$\tau(f)$ of $f$ is a section
of the induced vector bundle $f^{*}TN$
defined by
\be\tau(f)=
\sum_{i=1}^{n}\{\nabla^{f}_{e_i}df(e_i)
-df(\nabla_{e_i}e_i)\}\nonumber
\ee
for a local orthonormal frame $\{e_i\}$ on $M$, where  $\nabla^f$ and $\nabla$ denote
 the induced connection and the Levi-Civita connection of $M$, respectively.
If  $f$ is an isometric immersion,  then we have
\be\tau(f)=nH.\label{tau}
\ee

A smooth map $f$ is called
a {\it harmonic map} if 
it is a critical 
point of the energy functional
$$
E(f)=\int_{\Omega}|df|^2dv_g
$$
over every compact domain $\Omega$ of $M$, where $dv_g$ is the volume form of $M$. 
A smooth map $f$ is harmonic if and only if $\tau(f)$ vanishes identically on $M$.

\begin{definition}
{\rm For each smooth map $f$ of a compact domain $\Omega$ of $M$ into $N$, 
the $\lambda$-bienergy functional  is  defined by
\be
E_{2, \lambda}(f)=\int_{\Omega}|\tau(f)|^2dv_g+\lambda E(f).\nonumber
\ee
For simplicity, we denote $E_{2, 0}(f)$ by $E_2(f)$, which is  
called the bienergy functional.}
\end{definition}

Eliasson  \cite{eli}
proved that $E_{2, \lambda}$ satisfies 
Condition (C) of Palais-Smale if the dimension of the domain
is 2 or 3 and the target is non-positively curved. In general,  $E_{2, \lambda}$
does not satisfy Condition (C) (see \cite{le}).

%

\subsection{Ideal CR immersions as 
critical points of $\lambda$-bienergy functional}

Let $(M, HM, J_H, g)$ be a compact Riemannian almost CR manifold (with or without boundary) whose 
CR dimension is $h$, i.e.,
a compact smooth manifold equipped with a subbundle $HM$ of $TM$ of rank $2h$
together with a bundle isomorphism 
$J_H:HM\rightarrow HM$ such that $(J_H)^2=-I$, and a Riemannian metric $g$ such that
$g(X, Y)=g(J_HX, J_HY)$ for all $X$,$Y\in HM$.

An 
 immersion $f$ of   $(M, HM, J_H, g)$ into $\tilde M^m(4\epsilon)$
 is called a {\it CR immersion} if $J(df(X))=df(J_H(X))$ for any $X\in HM$.  
If $f$ is an isometric immersion, then $f(M)$ is a CR submanifold of $\tilde M^m(4\epsilon)$.
We denote  by
$\mathcal{ICR}(M, \tilde M^m(4\epsilon))$ the family of isometric  CR 
immersions of $M$ into $\tilde M^m(4\epsilon)$.
By Proposition \ref{pro3} and \ref{pro4}, we see that a
$\delta(n_1, \ldots, n_k)$-ideal  CR immersion of $M$ into $\tilde M^m(4\epsilon)$ is a 
stable critical point of $E_{2, \lambda}$
within  the class of $\mathcal{ICR}(M, \tilde M^m(4\epsilon))$.

\subsection{$\lambda$-biharmonic submanifolds and their extensions}
\begin{definition}[\cite{fel}]\label{e2}
{\rm A smooth map $f: M\rightarrow N$ is  called a
{\it $\lambda$-biharmonic map} if it is a critical point of the $\lambda$-bienergy
functional
with respect to all variations with compact support. 
If $f$ is a $\lambda$-biharmonic isometric immersion, then $M$  is called a {\it $\lambda$-biharmonic submanifold} in $N$. In the case of $\lambda=0$, we simply call it a 
{\it biharmonic submanifold}.}
\end{definition}

The Euler-Lagrange equation for $E_{2, \lambda}$ is given by (see \cite{ji2} and 
\cite[p.515]{fel})
\be
\tau_{2, \lambda}:=-\Delta_f(\tau(f))+{\rm trace}R^{N}(\tau(f),df)df-\lambda\tau(f)=0,\label{EL}\ee
where
$\Delta_f=-\sum_{i=1}^{n}(\nabla^{f}_{e_i}\nabla^{f}_{e_i}-\nabla^{f}_{\nabla_{e_i}e_i})$ and 
$R^{N}$ is 
the curvature tensor  of $N$, which is defined by $$R^{N}(X, Y)Z=[\nabla^N_X, \nabla^N_Y]Z-\nabla^N_{[X, Y]}Z$$ for the Levi-Civita connection $\nabla^N$ of $N$.
For simplicity, we denote  $\tau_{2, 0}(f)$ by $\tau_2(f)$.

By decomposing the left-hand side of 
(\ref{EL}) into  its tangential and normal components,  we have

\begin{proposition}[\cite{bal2}]\label{EL3}
Let $M$ be an $n$-dimensional 
submanifold of $\tilde R^m(\epsilon)$. Then $M$ is $\lambda$-biharmonic
if and only if 
\be
\begin{cases}
\Delta^{D}H+{\rm trace}\hskip2pt B(\cdot, A_H(\cdot))+(\lambda-\epsilon n)H=0,\\
4{\rm trace} A_{D_{(\cdot)}H}(\cdot)+n{\rm grad}(|H|^2)=0,\nonumber
\end{cases}
\ee
where $\Delta^D=-\sum_{i=1}^n\{D_{e_i}D_{e_i}-D_{\nabla_{e_i}e_i}\}$.
\end{proposition}

\begin{proposition}[\cite{fel}]\label{EL2}
Let $M$ be an $n$-dimensional 
submanifold of $\tilde M^m(4\epsilon)$ such that
$JH$ is tangent to $M$. Then $M$ is $\lambda$-biharmonic
if and only if 
\be
\begin{cases}
\Delta^{D}H+{\rm trace}\hskip2pt B(\cdot, A_H(\cdot))+\{\lambda-\epsilon(n+3)\}H=0,\\
4{\rm trace} A_{D_{(\cdot)}H}(\cdot)+n{\rm grad}(|H|^2)=0.\nonumber
\end{cases}
\ee
\end{proposition}

\begin{remark}
{\rm By Proposition \ref{EL2}, we see that
all  hypersurfaces with constant principal curvatures in $R^m(\epsilon)$ and  $\tilde M^m(4\epsilon)$
are $\{-|B|^2+\epsilon(m-1)\}$-biharmonic and $\{-|B|^2+2\epsilon(m+1)\}$-biharmonic, respectively.}
\end{remark}

It follows from (\ref{tau}) and (\ref{EL})  that any
minimal submanifold is $\lambda$-biharmonic. 
Thus, it is interesting to investigate non-minimal
$\lambda$-biharmonic submanifolds.
 
\begin{remark}{\rm Let $f: M\rightarrow \mathbb{R}^n$ be an isometric immersion.
We denote  the mean curvature vector field of $M$ by $H=(H_1, \ldots, H_n)$.
 Then, it follows from  (\ref{tau}) and (\ref{EL}) that $M$ is $\lambda$-biharmonic if and only if it satisfies 
 \be \Delta_{M}H_i=-\lambda H_i, \quad 1\leq i\leq n,\label{null}
 \ee
 where $\Delta_M$ is the Laplace operator acting on $C^{\infty}(M)$. 
Hence, the notion of biharmonic submanifolds in Definition \ref{e2}.
  is same as  one defined by B. Y. Chen (cf. \cite{chen9}).
 It was proved in \cite{chen8} 
   that 
   a submanifold $M$ satisfies (\ref{null}) 
 if and only if one of the following three cases occurs:
 
(1) $f$ satisfies $\Delta_Mf=-\lambda f$; 

(2) $f$ can be written as 
$f=f_0+f_1$, \quad $\Delta_Mf_0=0$, \quad$\Delta_Mf_1=-\lambda f_1$;

(3) $M$ is a biharmonic submanifold.

 An immersion described in (1) (resp. (2)) is said to be of {\it $1$-type}
 (resp. {\it null $2$-type}).
 An immersion $f: M\rightarrow \mathbb{R}^n$ is of $1$-type if and only if either $M$ is a minimal submanifold
 of  $\mathbb{R}^n$ or $M$ is a minimal submanifold of a hypersphere in $\mathbb{R}^n$
 (cf. \cite[Theorem 3.2]{chen7}). The classification of null $2$-type immersions 
 is not yet complete.

 }
\end{remark}
 
 There exist many non-minimal biharmonic submanifolds 
in a sphere or a complex projective space (see, for example, \cite{bal} and \cite{fel}).
On the other hand, the following conjecture proposed by Chen \cite{chen7} is still open.
\begin{conjecture}\label{conj}
{ Any biharmonic submanifold in  Euclidean space is minimal.}
\end{conjecture} 
Several 
  partial positive answers to this conjecture have been obtained (see \cite{chen9}).
  For example, Chen and Munteanu \cite{cm} proved that Conjecture \ref{conj}
  is true for  hypersurfaces which are $\delta(2)$-ideal or $\delta(3)$-ideal 
  in Euclidean space of 
  arbitrary dimension.


As an extension of the notion of biharmonic submanifolds, 
the following notion was introduced by Loubeau and Montaldo in \cite{lm}.
\begin{definition}\label{e3}
{\rm An isometric immersion $f: M\to N$
is called a   {\it $\lambda$-biminimal\/} 
if it is a critical point of
the $\lambda$-bienergy functional 
with respect to all
{\it normal variations} with compact support.
Here, a normal variation means a 
variation $f_t$ through $f=f_0$ 
such that the variational vector field 
$V=df_{t}/dt|_{t=0}$ is normal to $f(M)$. 
In this case, $M$ is called a {\it $\lambda$-biminimal submanifold} in $N$.  
In the case of  $\lambda=0$, we simply call it {\it biminimal} submanifold.}
\end{definition}
An isometric immersion $f$ is $\lambda$-biminimal if and only if 
\be
[\tau_{2, \lambda}(f)]^{\perp}=0,\nonumber
\ee
where $[\cdot]^{\perp}$ denotes the normal component of $[\cdot]$ (see \cite{lm}).
It is known that there exist ample examples of 
$\lambda$-biminimal submanifolds in real and complex space forms,
which are not $\lambda$-biharmonic
(see, for example, \cite{lm}, \cite{sasa7}, \cite{sa5} and \cite{sa6}).

 In \cite{sa}, the notion of  tangentially biharmonicity for submanifolds
   was introduced as follows:
  \begin{definition}\label{e4}
{\rm  Let $f:M\rightarrow N$ be an isometric immersion. Then $M$
 is called a {\it tangentially biharmonic submanifold} in $N$
 if it satisfies  
 \be
 [\tau_{2}(f)]^{\top}=0,\label{TBE}
 \ee where $[\cdot]^\top$ denotes the tangential part of $[\cdot]$. }
 \end{definition}



\begin{example}
{\rm Let $x:M^{n-1}\rightarrow \mathbb{R}^n$ be an isometric immersion.
The normal bundle $T^{\perp}M^{n-1}$ of $M^{n-1}$ is naturally immersed in $\mathbb{R}^n
\times\mathbb{R}^n=\mathbb{R}^{2n}$ by the immersion 
$f(\xi_x):=(x, \xi_x)$, which is  expressed as
\be
f(x, s)=(x, sV)\label{nbundle}
\ee
for the unit normal vector field $V$ along $x$.
We equip $T^{\perp}M^{n-1}$ with the metric induced by $f$. 
If we define the complex structure $J$ on $\mathbb{C}^n=\mathbb{R}^n\times\mathbb{R}^n$
by $J(X, Y):=(-Y, X)$, then 
$T^{\perp}M^{n-1}$ is a Lagrangian submanifold in $\mathbb{C}^n$
(see \cite[III.3.C]{hl}). 
It was proved  in \cite{sa} that  $T^{\perp}M^2$ is a tangentially biharmonic Lagrangian submanifold in
 $\mathbb{C}^3$
 if and only if
$M^2$ is either  minimal, a part of a round sphere or a part of a circular cylinder in $\mathbb{R}^3$.}
\end{example}

\begin{remark}
{\rm For any $\lambda\in\mathbb{R}$, we have $[\tau_{2, \lambda}(f)]^{\top}
=[\tau_{2}(f)]^{\top}$.}
\end{remark}

 \begin{remark}
 {\rm  By the first variation formula for $E_2$ obtained in \cite{ji2}, we see that 
 an isometric immersion $f:M\rightarrow N$ is tangentially biharmonic if and only if it 
 is a critical point of $E_{2}$ 
with respect to all
{\it tangential variations} with compact support.
Here, a tangential variation means a 
variation $f_t$ through $f=f_0$ 
such that the variational vector field 
$V=df_{t}/dt|_{t=0}$ is tangent to $f(M)$.}
 \end{remark}

 \begin{remark}
 {\rm As described by Hilbert \cite{hil}, 
 the stress-energy tensor associated to a variational problem is a symmetric $2$-covariant tensor which is conservative, namely, 
 divergence-free at critical points.
 The  stress-energy tensor $S_2$ 
for $E_{2}(f)$   was introduced by Jiang \cite{ji3} as follows:
\begin{align}
S_2(X, Y)=&\frac{1}{2}|\tau(f)|^2\<X, Y\>+\bigl<df, \nabla^{f}\tau(f)\bigr>\nonumber\\
&-\bigl<df(X), \nabla^{f}_Y\tau(f)\bigr>-
\bigl<df(Y), \nabla^{f}_X\tau(f)\bigr>.\nonumber
\end{align}
It satisfies ${\rm div}\hskip2pt S_2=-\<\tau_2(f), df\>$. Hence, 
an isometric 
immersion $f$ is tangentially biharmonic if and only if ${\rm div}\hskip2pt S_2=0$.
Caddeo et al. \cite{cad} called these submanifolds satisfying such a condition as
{\it biconservative submanifolds}, and moreover, 
 classified biconservative surfaces in 3-dimensional real space
forms.}
\end{remark}
\begin{remark}{\rm  
  Hasanis and Vlachos \cite{hav} classified 
  hypersurfaces in $\mathbb{R}^4$
  satisfying (\ref{TBE}). They called  such hypersurfaces as {\it $H$-hypersurfaces}.
  Afterwards, the biharmonic ones are picked out in the class.
  As a result, the non-existence of non-minimal
  biharmonic hypersurfaces in $\mathbb{R}^4$ was proved.}
  \end{remark}

\begin{remark}
{\rm It follows from  Definition \ref{e3} and \ref{e4}
 that a map $f$ is $\lambda$-biharmonic (resp. $\lambda$-biminimal) if and only if it 
is a critical point  of $E_{2}(f)$ for all 
variations (resp. normal variations) with compact support and {\it fixed energy}.  Here, 
$\lambda$ is the Lagrange multiplier.
}
\end{remark}

\subsection{Biharmonic ideal CR submanifolds}

For homogeneous real hypersurfaces in $\mathbb{C}P^m(4)$, namely, orbits
under some subgroups of the projective
unitary group $PU(m+1)$, we have

\begin{theorem}[\cite{iiu}]\label{clabiCP}
Let $M$ be a homogeneous hypersurface in $\mathbb{C}P^m(4)$.
Then, $M$ is non-minimal biharmonic if and only if it
is congruent to an open portion of one of the following real hypersurfaces{\rm :}

{\rm (1)} a tubular hypersurface  over $\mathbb{C}P^q(4)$ with radius  
\be
r=\cot^{-1}\biggl(\sqrt{\frac{m+2\pm\sqrt{(2q-m+1)^2+4(m+1)}}{2m-2q-1}}\biggr).\nonumber
\ee

{\rm (2)} a tubular hypersurface over the Pl\"{u}cker imbedding of the complex Grassmann manifold
$Gr_2(\mathbb{C}^5)\subset \mathbb{C}P^9(4)$ with  radius $r$, where $0<r<\pi/4$ and 
$t=\cot{r}$ is a unique solution of the equation
\be 41t^6+43t^4+41t^2-15=0.\nonumber\ee
In this case, $r=1.0917\cdots$.

{\rm (3)} a tubular hypersurface over the canonical imbedding of the 
Hermitian symmetric space $SO(10)/U(5)\subset \mathbb{C}P^{15}(4)$ 
with radius $r$, where $0<r<\pi/4$ and $t=\cot{r}$
is a unique solution of the equation
\be
13t^6-107t^4+43t^2-9=0.\nonumber\ee
In this case, $r=0.343448\cdots$. 
\end{theorem}

For details on the canonical imbedding of a compact Hermitian symmetric space into
$\mathbb{C}P^m(4)$, we refer the reader to  Section 4 of \cite{nt}. 
\begin{remark}
{\rm Let $M$ be a real hypersurface in $\mathbb{C}P^m(4)$. Kimura \cite{kim} proved that
 $M$ is a Hopf 
hypersurface with constant principal curvatures
if and only if it is homogeneous.}
\end{remark}

Combining Theorem \ref{claidCP}, Proposition \ref{EL2} and Theorem \ref{clabiCP}, we obtain  
\begin{corollary}
Let $M$ be a $\delta_k(2)$-ideal non-minimal biharmonic Hopf hypersurface in $\mathbb{C}P^m(4)$.
Then, $m$ is odd and 
$M$ is an open portion of  a tubular hypersurface  over $\mathbb{C}P^{(m-1)/2}(4)$ with radius  
\be
r=\cot^{-1}\biggl(\sqrt{\frac{m+2\pm 2\sqrt{m+1}}{m}}\biggr).\nonumber
\ee
\end{corollary}

\begin{example}{\rm 
On each CR submanifold described in Theorem \ref{cla},  
there exists an orthonormal frame $\{e_1, \ldots, e_{2m}\}$ such that $e_{2r}=Je_{2r-1}$ 
for $r\in\{1, \ldots, n\}$, $JH\parallel e_{2n+1}\in
\mathcal{H}^{\perp}$ and the second fundamental form $B$ takes the following form:
\begin{align}
&B(e_{2r-1}, e_{2r-1})=\sqrt{\frac{k}{2n-k}}Je_{2n+1}+\phi_r\xi_r,\nonumber\\
&B(e_{2r}, e_{2r})=\sqrt{\frac{k}{2n-k}}Je_{2n+1}-\phi_r\xi_r,\nonumber\\
&B(e_{2r-1}, e_{2r})=\phi_rJ\xi_r,\nonumber\\
&B(e_{2n+1},e_{2n+1})=\frac{2n}{\sqrt{k(2n-k)}}Je_{2n+1},\nonumber\\
&B(u_i, u_j)=h(u_i, e_{2n+1})=0\quad (i\ne j),\nonumber
\end{align}
where $\phi_r$ are functions, 
$\xi_r\in\nu$ and  $u_j\in L_j$ (see Lemma 7 of \cite{sa3}).
Here,  $\nu$ denotes an orthogonal complement
 of $J\mathcal{H}^{\perp}$ in $T^{\perp}M$.
Therefore, by using Proposition \ref{EL2}, we find that 
all ideal CR submanifolds given in Theorem \ref{cla}
are  non-minimal $\lambda$-biharmonic submanifolds with
\be
\lambda=-\frac{2n(2n+k^2)}{k(2n-k)}-2n-4\hskip5pt (\ne 0).\nonumber
\ee}
\end{example}

The following problem seems interesting.
\begin{problem}
Classify $\delta(n_1, \ldots, n_k)$-ideal proper
CR submanifolds   in $\mathbb{C}P^m(4)$  which are non-minimal biharmonic.
\end{problem}

\begin{example}
{\rm The standard product $S^{2r+1}(1/\sqrt{2})\times S^{2s+1}(1/\sqrt{2})$  in $S^{2(r+s)+3}(1)$
 is a   biharmonic contact CR hypersurface (see \cite[Example 5.1]{yano}
 and \cite[p.92]{bal}).
 Its
 principal curvatures  are  $\{1, -1\}$ with multiplicities $\{2r+1, 2s+1\}$. We may assume that $r\geq s$.
 By Theorem \ref{ingere}, we see that
 the biharmonic hypersurface $S^{2r+1}(1/\sqrt{2})\times S^{2s+1}(1/\sqrt{2})$ 
 is minimal and $\delta_{2s+1}(2)$-ideal if $r=s$; otherwise it is
  non-minimal and $\delta(4s+3)$-ideal. }
  \end{example}

 Incidentally, the following conjectures proposed in \cite{bal2} remains open.
 \begin{conjecture}
 The only non-minimal biharmonic hypersurfaces in $S^{m+1}$ are the open
 parts of hyperspheres $S^{m}(1/\sqrt{2})$ or of the standard products $S^{m_1}(1/\sqrt{2})\times
 S^{m_2}(1/\sqrt{2})$, where $m_1+m_2=m$ and $m_1\ne m_2$.
 \end{conjecture}

\begin{conjecture}
Any non-minimal biharmonic submanifold in $S^n(1)$ has constant mean curvature.
\end{conjecture}







 \vskip20pt
 
 
 
 


 \end{document}